\documentclass[12pt]{article}
\usepackage{latexsym,amsfonts,amssymb}
\usepackage[dvips]{graphics}
\begin{document}
\footnotetext{E-mail address: weigenyan@263.net(Weigen Yan),
mayeh@math.sinica.edu.tw(Y.-N Yeh),
fjzhang@jingxian.xmu.edu.cn(Fuji Zhang)\\ $^a$ Partially supported
by FMSTF(2004J024) and NSFF(E0540007). $^b$ Partially supported by
NSC94-2115-M001-017. $^c$Partially supported by NSFC (10371102).}
\begin{center}
{\Large\bf Graphical condensation of plane graphs: a combinatorial
approach}
\end{center}
\vskip0.2cm
\begin{center}
{\small \it Weigen Yan$^{1,2,a}$,\ \ \ \ Yeong-Nan Yeh$^{2,b}$,\ \ \ \ Fuji Zhang$^{3,c}$\\
\small
(1. School of Sciences, Jimei University, Xiamen 361021, China)\\
\small (2. Institute of Mathematics, Academia Sinica,
Taipei 11529, Taiwan)\\
 \small (3. Department of Mathematics, Xiamen University, Xiamen 361005, China) }\\
\end{center}
\begin{center}
\begin{minipage}{135mm}
\begin{center}
{\bf Abstract}
\end{center}
\ \ \ \ \ {\small The method of graphical vertex-condensation for
enumerating perfect matchings of plane bipartite graph was found
by Propp (Theoret. Comput. Sci. 303(2003), 267-301), and was
generalized by Kuo (Theoret. Comput. Sci. 319 (2004), 29-57) and
Yan and Zhang (J. Combin. Theory Ser. A, 110(2005), 113-125). In
this paper, by a purely combinatorial method some explicit
identities on graphical vertex-condensation for enumerating
perfect matchings of plane graphs (which do not need to be
bipartite) are obtained. As applications of our results, some
results on graphical edge-condensation for enumerating perfect
matchings are proved, and
%a new proof of Stanley's multivariate
%version of the Aztec diamond theorem is given,
we count the sum of weights of perfect matchings of weighted Aztec
diamond. }
\end{minipage}
\end{center}
\baselineskip 0.1in\pagestyle{plain} \setcounter{page}{1}
\section* {\bf 1. Introduction}
\label{sec:intro}
\label{sec:subdivisions}

\par \  \ \ Throughout this paper, we suppose that
$G=(V(G),E(G))$ is a simple graph with the vertex set
$V(G)=\{v_1,v_2,\cdots,v_n\}$ and the edge set
$E(G)=\{e_1,e_2,\cdots, e_m\}$, if not specified. A perfect
matching of $G$ is a set of independent edges of $G$ covering all
vertices of $G$. Denote the set of perfect matchings of $G$ by
$\mathcal {M}(G)$ and the number of perfect matchings of $G$ by
$M(G)$. If $G$ is a weighted graph, the weight of a perfect
matching $P$ of $G$ is defined to be the product of weights of
edges in $P$. We also denote the sum of weights of perfect
matchings of $G$ by $M(G)$. Let $A=\{a_1,a_2,\ldots,a_s\}$ (resp.
$E_1=\{e_{i_1}, e_{i_2},\ldots, e_{i_t}\}$) be a subset of the
vertex set $V(G)$ (resp. a subset of the edge set $E(G)$). By
$G-A$ or $G-a_1-a_2-\ldots-a_s$ (resp. $G-E_1$ or
$G-e_{i_1}-e_{i_2}-\ldots-e_{i_t}$) we denote the induced subgraph
of $G$ by deleting all vertices in $A$ and the incident edges from
$G$ (resp. by deleting all edges in $E_1$).\\

\par By the method of graphical condensation for
enumerating perfect matchings of plane bipartite graphs,
Propp \cite{Propp03} obtained the following result.\\

{\bf Proposition 1.1} (Propp \cite{Propp03})\ Let $G=(U,V)$ be a
plane bipartite graph in which $|U|=|V|$. Let vertices $a, b, c$
and $d$ form a $4-$cycle face in $G$, $a, c\in U$, and $b, d\in
V$. Then
$$M(G)M(G-\{a,b,c,d\})=M(G-\{a,b\})M(G-\{c,d\})+M(G-\{a,d\})M(G-\{b,c\}).$$

\par By a combinatorial method, Kuo \cite{Kuo04} generalized Propp's result above as follows.\\

{\bf Proposition 1.2} (Kuo \cite{Kuo04})\ Let $G=(U,V)$ be a plane
bipartite graph in which $|U|=|V|$. Let vertices $a, b, c,$ and
$d$ appear
in a cyclic order on a face of $G$.\\
{\bf (1)}\ \ If $a, c\in U$, and $b, d\in V$, then
$$M(G)M(G-\{a,b,c,d\})=M(G-\{a,b\})M(G-\{c,d\})+M(G-\{a,d\})M(G-\{b,c\}).$$
{\bf (2)}\ \ If $a, b\in U$, and $c, d\in V$, then
$$M(G-\{a,d\})M(G-\{b,c\})=M(G)M(G-\{a,b,c,d\})+M(G-\{a,c\})M(G-\{b,d\}).$$

\par By Ciucu's Matching Factorization Theorem in \cite{Ciucu97}, Yan and Zhang
\cite{YZ05}
obtained a more general result than Kuo's as follows.\\

{\bf Proposition 1.3} (Yan and Zhang \cite{YZ05})\ \ Let $G=(U,V)$
be a plane weighted bipartite graph in which $|U|=|V|=n$. Let
vertices $a_1, b_1, a_2, b_2, \ldots, a_k, b_k\ (2\leq k\leq n)$
appear in a cyclic order on a face of $G$, and let $A_1=\{a_i\ |\
a_i\in U,\ 1\leq i\leq k\}$, $A_2=\{a_i\ |\ a_i\in V, \ 1\leq
i\leq k\}$, $B_1=\{b_i\ |\ b_i\in V, \ 1\leq i\leq k \}$ and
$B_2=\{b_i\ |\ b_i\in U, \ 1\leq i\leq k\}$. If $|A_1\cup
B_2|=|A_2\cup B_1|=k$, then
$$2^kM(G-A_1-B_1)M(G-A_2-B_2)=\sum_{(X,Y)\subseteq (A_1\cup B_2)
\times (A_2\cup B_1),|X|=|Y|} M(G-X-Y)M(G-\overline X-\overline
Y),$$ where the sum ranges over all subsets $(X,Y)$ of $(A_1\cup
B_2) \times (A_2\cup B_1)$ such that $|X|=|Y|$, and $X\subseteq
(A_1\cup B_2), Y\subseteq (A_2\cup B_1), \overline X=(A_1\cup B_2)
\backslash X, \overline
Y=(A_2\cup B_1)\backslash Y$.\\

\par The results above hold under the condition that the plane graph
considered is bipartite. For the case in which the plane graph
does not need to be bipartite, in an email sent to ``Domino Forum"
Propp wrote that Kenyon recently told him about an identity of
Pfaff's that, in combination with Kasteleyn's Pfaffian method (see
\cite{Kas63,Kas67}), implies the
following combinatorial assertion:\\

{\bf Proposition 1.4}\ \ Let G be a plane graph with four vertices
a,b,c,d (in the cyclic order) adjacent to a single face. Then

$$
M(G)M(G-\{a,b,c,d\})+M(G-\{a,c\})M(G-\{b,d\})
$$
$$
=M(G-\{a,b\})M(G-\{c,d\})+M(G-\{a,d\})M(G-\{b,c\}).\eqno{(1)}$$
Propp also hoped to find a combinatorial proof of $(1)$. Kuo told
a result similar to Proposition 1.4 in ``Domino Forum".
%Later, in
%``Domino Forum" Krattenthaler claimed that Kuo's result and (1)
%were all special instances of a family of Pfaffian identities,
%apparently discovered for the first time by Wenzel [16], Dress and
%Wenzel [5] (which was also proved by a new method in Knuth [10]),
%and further generalized by Hamel [7]. Particularly, Krattenthaler
%said that we could generate many more identities similar to (1).
But it seems that the explicit results (including the identity
(1)) have not been published. Furthermore, it seems that nobody
has published a purely combinatorial proof of $(1)$.

In the next section, inspired by an interesting lemma in Ciucu
\cite{Ciucu97} and some Pfaffian identities (see
\cite{DW95,Ham01,Knuth96,Wenzel93}), we find a purely
combinatorial method to obtain some explicit identities concerning
the enumeration of perfect matchings of plane graphs, which do not
need to be bipartite. Our results imply Propositions 1.2 and 1.4.
On the other hand, an obvious observation in the identities in
Propositions $1.1-1.4$ is that the graphs related in these
identities are either $G$ or the induced subgraphs of $G$ by
deleting some vertices. For the sake of convenience, we call these
procedures for enumerating perfect matchings ``graphical
vertex-condensation" in place of ``graphical condensation", the
term used by Kuo \cite{Kuo04}. In other words, we regard Kuo's
``graphical condensation" as ``condensing vertices of bipartite
graphs ". Based on this, it is natural to ask whether we can
condense edges of $G$ or both of edges and vertices. The theorems
and corollaries in Section 3 answer this question in the
affirmative. We call these results ``graphical edge-condensation"
for enumerating perfect matchings of plane graphs. In Section 4,
we obtain a new proof of Stanley's multivariate version of the
Aztec diamond theorem.

\section* { \bf 2. Graphical vertex-condensation
%for enumerating perfect matchings of plane graphs
} \ \ \ \ We say a plane graph $G$ is symmetric if it is invariant
under the reflection across some straight line $\ell$ (say
symmetry axis). Figure 1(a) shows an example of a symmetric plane
graph. A weighted symmetric graph is a symmetric graph equipped
with weight on every edge of $G$ that is constant on the orbits of
the reflection. The width of a symmetric graph $G$, denoted by
$\omega (G)$, is defined to be half the number of vertices of $G$
lying on the symmetric axis. Clearly, if $\omega (G)$ is not an
integer then $M(G)=0$. Hence we suppose that there are even number
of vertices of $G$ lying on the symmetry axis.
%%%%%%%%%%%%%%%%%%%%%%%%%%%%%%%%%%%%%%%%%%
%%%%%%%%%%%%% Figure 1
%%%%%%%%%%%%%%%%%%%%%%%%%%%%%%%%%%%%%%%%%%
\begin{figure}[htbp]
  \centering
  \includegraphics{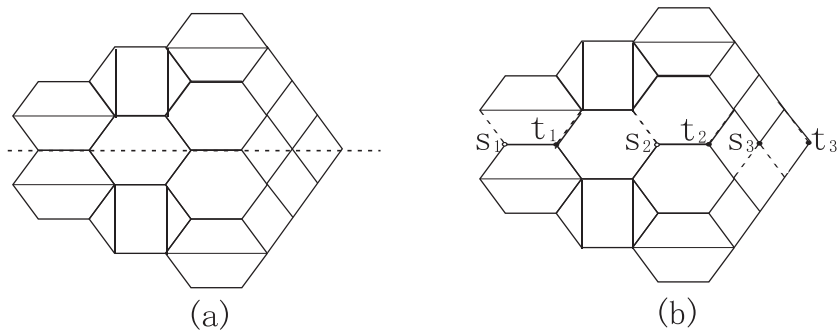}
  \caption{\ (a)\ A symmetric graph $G$. \ (b)\ A
reduced subgraph of symmetric graph $G$.}
\end{figure}

\par Let $G$ be a plane weighted symmetric graph with symmetry axis $\ell$, which we
consider to be horizontal. Let $s_1, t_1, s_2, t_2, \cdots, s_k,
t_k$ be the vertices lying on $\ell$ as they occur from left to
right. A reduced subgraph of $G$ is a graph obtained from $G$ by
deleting at each vertex $s_i$ either all incident edges above
$\ell$ or all incident edges below $\ell$. Figure 1(b) shows a
reduced subgraph of the graph presented in Figure 1(a) (the
deleted edges of the original graph are represented by dotted
lines). Obviously, there exist exactly $2^k$ reduced subgraphs of
$G$. Now, we can introduce a lemma found by Ciucu \cite{Ciucu97}
and proved by a purely combinatorial method, which plays a key
role in the proof of one of our main theorems.\\

{\bf Lemma 2.1} (Ciucu \cite{Ciucu97})\ \ Let $G$ be a plane
weighted symmetric graph and there exist $2k$ vertices lying on
the symmetry axis. Then all $2^k$ reduced subgraphs of $G$ have
the same sum of weights of perfect matchings.\\

\par Now we are in the position to prove one of our main results.\\

{\bf Theorem 2.2}\ \ Let $G$ be a plane weighted graph with $2n$
vertices. Let vertices $a_1, b_1, a_2, b_2, \ldots, a_k, b_k\
(2\leq k\leq n)$ appear in a cyclic order on a face of $G$, and
let $A=\{a_1, a_2, \cdots, a_k\}$, $B=\{b_1, b_2, \cdots, b_k\}$.
Then, for any $j=1, 2, \cdots, k$, we have
$$\sum_{Y\subseteq B,\ |Y|\ \mbox {is odd}} M(G-a_j-Y)M(G-A\backslash \{a_j\}-\overline
Y)$$
$$=\sum_{W\subseteq B,\ |W|\ \mbox{is even}}M(G-W)M(G-A-\overline W).\eqno{(2)}$$
where the first sum ranges over all odd subsets $Y$ of $B$ and the
second sum ranges over all even subsets $W$ of $B$, $\overline
Y=B\backslash Y$ and $\overline W=B\backslash W$.\\

{\bf Proof}\ \ Since $G$ is a plane graph, for an arbitrary face
$F$ of $G$ there exists a planar embedding of $G$ such that the
face $F$ is the unbounded one. Hence we may assume that vertices
$a_1, b_1, a_2, b_2, \ldots, a_k, b_k$ appear in a cyclic order on
the unbounded face of $G$. Take two copies of the weighted graph
$G$, denoted by $G_1=(V(G_1),E(G_1))$ with the vertex set
$V(G_1)=\{v_i^{(1)}|\ 1\leq i\leq 2n\}$, and $G_2=(V(G_2),E(G_2))$
with the vertex set $V(G_2)=\{v_i^{(2)}|\ 1\leq i\leq 2n\}$,
respectively, and leave weights of all edges unchanged. Hence
$a_1^{(1)}, b_1^{(1)}, a_2^{(1)}, b_2^{(1)}, \ldots, a_k^{(1)},
b_k^{(1)}$ appear in a cyclic order on the unbounded face of $G_1$
and $a_1^{(2)}, b_1^{(2)}, a_2^{(2)}, b_2^{(2)}, \ldots,
a_k^{(2)}, b_k^{(2)}$ appear in a cyclic order on the unbounded
face of $G_2$. Construct a new plane weighted graph with $4n+2k$
vertices, denoted by $\widetilde G=(V(\widetilde G), E(\widetilde
G ))$, such that $V(\widetilde G)=V(G_1)\cup V(G_2)\cup W$,
$E(\widetilde G)=E(G_1)\cup E(G_2)\cup
\{a_i^{(1)}s_i,a_i^{(2)}s_i,b_i^{(1)}t_i,b_i^{(2)}t_i|\ 1\leq
i\leq k \}$, where $W=\{s_1,t_1,s_2,t_2,\ldots,s_k,t_k\}$. Let the
weight of every edge in
$\{a_i^{(1)}s_i,a_i^{(2)}s_i,b_i^{(1)}t_i,b_i^{(2)}t_i|\ 1\leq
i\leq k \}$ in $\widetilde G$ be 1 and leave all other weights
unchanged. The resulting weighted graph is $\widetilde G$. Figure
2(a) and (b) show this procedure constructing the new weighted
graph $\widetilde G$ from the weighted graph $G$. Obviously,
$\widetilde G$ is a plane weighted graph. Furthermore, by the
definition of the symmetric graph, $\widetilde G$ can be regarded
as a symmetric weighted plane graph with symmetry axis $\ell$,
which contains $2k$ vertices lying on $\ell$.
%%%%%%%%%%%%%%%%%%%%%%%%%%%%%%%%%%%%%%%%%%
%%%%%%%%%%%%% Figure 2
%%%%%%%%%%%%%%%%%%%%%%%%%%%%%%%%%%%%%%%%%%
\begin{figure}[htbp]
  \centering
  \scalebox{0.8}{\includegraphics{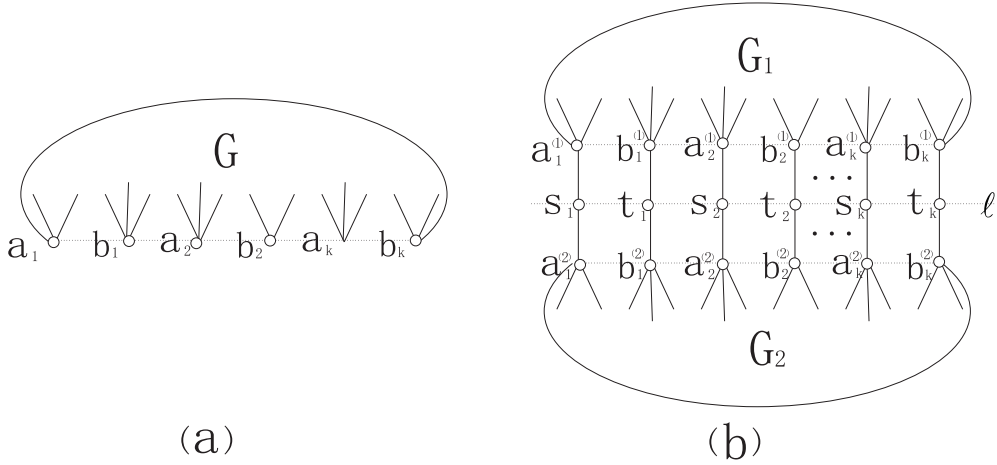}}
  \caption{\ (a)\ The graph $G$. \ (b)\ The graph $\widetilde G$.}
\end{figure}

Now, we consider the following $k+1$ reduced subgraphs of
$\widetilde G$, denoted by $G^{(0)}, G^{(1)}, \cdots, G^{(k)}$,
respectively, where $G^{(i)}=\widetilde G-E_i$,
$E_0=\{s_pa_p^{(1)}|\ p=1,2,\cdots,k\}, E_i=\{s_pa_p^{(1)}|\
p=1,2,\cdots,i-1,i+1,\cdots, k\}\cup \{s_ia_i^{(2)}\}$ for
$i=1,2,\cdots,k$. Hence, by Lemma 2.1, we have
$$M(G^{(0)})=M(G^{(1)})=\cdots=M(G^{(k)}).\eqno{(3)}$$
We partition the set $\mathcal M(G^{(0)})$ of perfect matchings of
$G^{(0)}$ such that
$$\mathcal M(G^{(0)})=\mathcal M_0\cup \mathcal M_1\cup \ldots \cup
\mathcal M_{[\frac{k}{2}]},$$ where $\mathcal M_i$ denotes the set
of perfect matchings of $G^{(0)}$ containing exactly $2i$ edges in
subset $\{t_jb_j^{(1)}|\ 1\leq j\leq k\}$ of $E(G^{(0)})$. It is
obvious that, for any $i$ ($0\leq i\leq [\frac{k}{2}]$), after
removing the forced edges we have
$$|\mathcal M_i|=\sum_{Y\subseteq B,\ |Y|=2i} M(G-Y)M(G-A-\overline Y),$$ where the sum
ranges over all subsets $Y$ of $B$ such that $|Y|=2i$.
%Particularly, $|\mathcal M_0|=M(G)M(G-A-B)$.
Hence we have
$$M(G^{(0)})=|\mathcal M(G^{(0)})|=\sum_{i=0}^{[\frac{k}{2}]}|\mathcal M_i|
=\sum_{Y\subseteq B,\ |Y|\ \mbox{is even }}M(G-Y)M(G-A-\overline
Y),\eqno{(4)}$$ where the second sum
ranges over all even subsets of $B$.\\

\par Similarly, for any $j=1,2,,\cdots,k$, we can prove that
$$M(G^{(j)})=
\sum_{Y\subseteq B,\ |Y|\ \mbox {is odd}}
M(G-a_j-Y)M(G-A\backslash \{a_j\}-\overline Y),\eqno{(5)}$$ where
the sum ranges over all odd subsets of $B$.\\

\par The theorem thus follows from (3)$-$(5).\ \ \ $\blacksquare$\\

{\bf Remark 1}\ \ Note that Ciucu \cite{Ciucu97} used a purely
combinatorial method to prove Lemma 2.1. Hence, by the procedure
proving Theorem
2.2, our method to prove Theorem 2.2 is also combinatorial.\\

{\bf Remark 2}\ \ Proposition 1.4 is the special case of Theorem
2.2 in which $k=2$.\\

\par The following corollary, which has a simpler form than
that in Corollary 2.3 in Yan and Zhang \cite{YZ05}, is
the special instance of Theorem 2.2.\\

{\bf Corollary 2.3} (Yan and Zhang \cite{YZ05})\ \ Let $G=(U,V)$
be a plane weighted bipartite graph in which $U=\{u_i|1\leq i\leq
n\}$ and $V=\{v_i|1\leq i\leq n\}$. Let vertices $a_1, b_1, a_2,
b_2, \ldots, a_k, b_k$ appear in a cyclic order on a face of $G$.
If $A=\{a_i|1\leq i\leq k\}\subseteq U$, and $B=\{b_i|1\leq i\leq
k\}\subseteq V$, then
$$M(G)M(G-A-B)=\sum_{i=1}^n
M(G-a_j-b_i)M[G-(A\cup B)\backslash \{a_j,b_i\}]\eqno{(6)}$$ for
any $j=1,2,
\ldots,k$. \\

{\bf Proof}\ \ Note that $G=(U,V)$ is a bipartite graph, and
$A=\{a_i|1\leq i\leq k\}\subseteq U$ and $B=\{b_i|1\leq i\leq
k\}\subseteq V$. Hence, in the formula $(2)$ in Theorem 2.2 if
$|Y|$ is an odd integer more than $1$ we have $M(G-a_j-Y)=0$.
Similarly, in the formula $(2)$ in Theorem 2.2 if $|W|\neq 0$ we
have $M(G-W)=0$. Thus it is not difficult to see that $(6)$ is
immediate from $(2)$. \ \ \ \ $\blacksquare$\\

If we set $k=3$ in Corollary 2.3, then we have the following
formula:
$$M(G)M(G-a_1-a_2-a_3-b_1-b_2-b_3)=M(G-a_1-b_1)M(G-a_2-a_3-b_2-b_3)+$$$$M(G-a_1-b_2)M(G-a_2-a_3-b_1-b_3)
+M(G-a_1-b_3)M(G-a_2-a_3-b_1-b_2),\eqno{(7)}$$

\par {\bf Remark 3}\ \ Similarly, we can obtain the identities in Corollaries 2.5 and 2.6 in Yan and Zhang
\cite{YZ05} from Theorem 2.2.

\section* { \bf 3. Graphical edge-condensation
%for enumerating perfect matchings of plane graphs
}

\ \ \ \ Let $G=(V(G),E(G))$ be a weighted graph and $e=ab$ an edge
of $G$. Define a new weighted graph $G'=(V(G'),E(G'))$ from $G$ as
follows. Delete the edge $e=ab$ from $G$ and add three edges $aa',
a'b', b'b$ with the weights $\sqrt {\omega_e}, 1$ and $\sqrt
{\omega_e}$, where $\omega_e$ denotes the weight of edge $e$. The
resulting weighted graph is $G'$. Hence $V(G')=\{a',b'\}\cup V(G)$
and $E(G')=\{aa',a'b',b'b\}\cup E(G)\backslash \{e\}$. Figure 3
(a) and (b) illustrate this procedure.
%%%%%%%%%%%%%%%%%%%%%%%%%%%%%%%%%%%%%%%%%%
%%%%%%%%%%%%% Figure 3
%%%%%%%%%%%%%%%%%%%%%%%%%%%%%%%%%%%%%%%%%%
\begin{figure}[htbp]
  \centering
  \scalebox{0.8}{\includegraphics{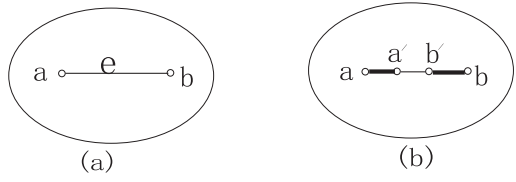}}
  \caption{\ (a)\ The weighted graph $G$ in
Lemma 3.1. \ (b)\ The weighted graph $G'$ obtained from $G$ in
Lemma 3.1.}
\end{figure}

{\bf Lemma 3.1} (Ciucu \cite{Ciucu96})\ \ Let $G$ be a weighted
graph and $e=ab$ an edge of $G$, and let $G'$ be the weighted
graph defined above. Then
$$M(G)=M(G').$$

\par In order to state our main results, we need to introduce some
notation. We use ${\bf [k]}$ to denote the set $\{1,2,\ldots,k\}$.
Let $G$ be a graph, and let $e_1=a_1b_1, e_2=a_2b_2, \ldots,
e_k=a_kb_k$ ($2\leq k\leq n$) be $k$ independent edges (a matching
of $G$ with $k$ edges) in $G$, and $X\subseteq A=\{a_i|1\leq i\leq
k\},\ Y\subseteq B=\{b_i|1\leq i\leq k\}$. Define: $I_X=\{i|a_i\in
X\}$, $I_Y=\{i|b_i\in Y\}$. Let $I$ be a subset of ${\bf [k]}$ and
$\bar I={\bf [k]}\backslash I$. Define: $E_I=\{e_i|\ i\in I\}$,
$A_I=\{a_i|\ i\in I\}$, $B_I=\{b_i|\ i\in I\}$. Let $I_1\subseteq
{\bf [k]}$ and
$I_2\subseteq {\bf [k]}$. Define: $I_1-I_2=I_1\backslash
(I_1\cap I_2), I_1\triangle I_2=(I_1-I_2)\cup (I_2-I_1)$.\\

{\bf Theorem 3.2}\ \ Suppose $G$ is a plane weighted graph with
even number of vertices and the weight of every edge $e$ in $G$ is
denoted by $\omega_e$. Let $e_1=a_1b_1, e_2=a_2b_2, \ldots,
e_k=a_kb_k$ ($k\geq 2$) be k independent edges in the boundary of
a face $f$ of $G$, and let vertices $a_1, b_1, a_2, b_2, \ldots,
a_k, b_k$ appear in a cyclic order on $f$, and let $A=\{a_i|\
i=1,2,\ldots,k\}$, $B=\{b_i|\ i=1,2,\ldots,k\}$ and $E=\{e_i|\
i=1,2,\ldots,k\}$. Then, for any $j=1,2,\ldots,k$,
\[\sum_{{W\subseteq B}\atop {|W|\ \mbox{\small is even}}}\left(\prod_{e\in
E_{I_W}}\omega_e\right)M(G-A_{I_W})M(G-E_{\overline{I_W}}-B_{I_W})=
\sum_{{Y\subseteq B}\atop {|Y|\ \mbox{\small is odd}}}
\left(\prod_{e\in E_{\{j\}\triangle I_Y}}\omega_e\right)\times\]
\[\left\{M(G-E_{I_Y\cap\{j\}}-B_{\{j\}-I_Y}-A_{I_Y-\{j\}})M(G-E_{\overline {I_Y}\cap \overline {\{j\}}}
-B_{\overline {\{j\}}-\overline {I_Y}}-A_{\overline
{I_Y}-\overline {\{j\}}})\right\}\ \ \ \ (8)\] where the first
product is over all edges in $E_{I_W}$, the second product is over
all edges in $E_{\{j\}\triangle I_Y}$, the first sum ranges over
all even subsets of $B$, and the
second sum ranges over all odd subsets of $B$.\\

{\bf Proof}\ \ Let $G'$ be the graph obtained from $G$ by deleting
$k$ edges $e_1, e_2, \ldots, e_k$ and adding $3k$ edges $a_ia_i',
a_i'b_i', b_i'b_i$ with the weights $\sqrt {\omega_{e_i}}, 1,
\sqrt {\omega_{e_i}}$ for $i=1, 2, \ldots, k$, and leaving all
other weights unchanged. Hence, the vertex set of $G'$, denoted by
$V(G')$, is $\{a_i',b_i'|1\leq i\leq k\}\cup V(G)$, and the edge
set of $G'$, denoted by $E(G')$, is
$\{a_ia_i',a_i'b_i',b_i'b_i|i=1,2,\ldots,k\}\cup E(G)\backslash
\{e_i|1\leq i\leq k\}$, where $V(G)$ and $E(G)$ are the vertex set
and the edge set of $G$, respectively. For the sake of
convenience, denote the edge $a_i'b_i'$ by $e_i'=a_i'b_i'$ for
$i=1, 2, \ldots, k$. Figure 4 (a) and (b) show this procedure.
%%%%%%%%%%%%%%%%%%%%%%%%%%%%%%%%%%%%%%%%%%
%%%%%%%%%%%%% Figure 4
%%%%%%%%%%%%%%%%%%%%%%%%%%%%%%%%%%%%%%%%%%
\begin{figure}[htbp]
  \centering
  \scalebox{0.8}{\includegraphics{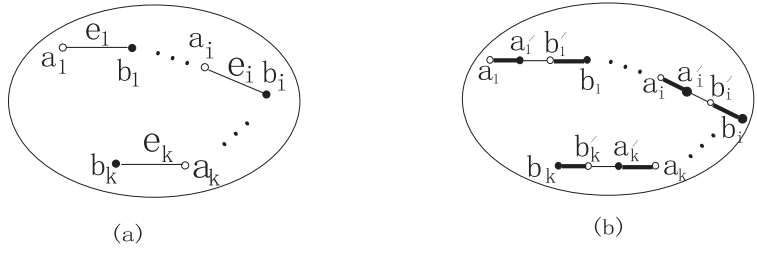}}
  \caption{\ (a)\ The weighted graph $G$ in
the proof of Theorem 3.2. \ (b)\ The weighted graph $G'$ obtained
from $G$ in the proof of Theorem 3.2.}
\end{figure}

Obviously, by the definition of $G'$, $G'$ is a plane weighted
graph with even number of vertices. Furthermore, vertices
$a_1',b_1',a_2',b_2',\ldots,a_k',b_k'$ appear in a cyclic order on
a face of $G'$.  Let $A'=\{a_i'|\ 1\leq i\leq k\}$ and
$B'=\{b_i'|\ 1\leq i\leq k\}$. By Theorem 2.2, we have
$$\sum_{{W'\subseteq B'}\atop {|W'|\ \mbox{is
even}}}M(G'-W')M(G'-A'-\overline {W'})$$$$=\sum_{{Y'\subseteq
B'}\atop {|Y'|\ \mbox{is odd}}} M(G'-a_j'-Y') M(G'-A'\backslash
\{a_j'\}-\overline {Y'}) \eqno{(9)}$$ for any $j=1,2\ldots,k$,
where the first sum ranges over all even subsets $W'$ of $B'$ and
the second sum is over all odd subsets $Y'$ of $B'$, and
$\overline {Y'}=B'\backslash Y', \overline {W'}=B'\backslash
{Y'}$.

\par Let $Y'$ be an odd subset of $B'$. By our notation defined
above, $I_{Y'}=\{i|\ b_i'\in Y'\}$. Let $Y=\{b_i|\ i\in I_{Y'}\}$.
Hence $I_Y=I_{Y'}$. Note that
$$M(G'-a_j'-Y')=M(G'-\{a_i',b_i'|\ i\in I_Y\cap \{j\}\}-\{a_i'|\ i\in \{j\}-I_Y\}-
\{b_i'|\ i\in I_Y-\{j\}\}).$$  By Lemma 3.1, after removing the
forced edges we have
$$M(G'-a_j'-Y')=\left(\prod_{e\in E_{(\{j\}-I_Y)\cup (I_Y-\{j\})}}\sqrt {\omega_e}\right)
M(G-E_{I_Y\cap\{j\}}-B_{\{j\}-I_Y}-A_{I_Y-\{j\}})$$
$$=\left(\prod_{e\in E_{\{j\}\triangle I_Y}}\sqrt {\omega_e}\right)
M(G-E_{I_Y\cap\{j\}}-B_{\{j\}-I_Y}-A_{I_Y-\{j\}}). \eqno{(10)}$$
Similarly, we have
$$M(G'-A'\backslash\{a_j'\}-\overline{Y'})=M(G'-\{a_i',b_i'|\ i\in \overline {\{j\}}\cap
\overline{I_Y}\} -\{a_i'|\ i\in
\overline{\{j\}}-\overline{I_Y}\}-\{b_i'|\ i\in
\overline{I_Y}-\overline{\{j\}}\})$$
$$=\left(\prod_{e\in E_{\overline{\{j\}}\triangle \overline{I_Y}}}
\sqrt {\omega_e}\right) M(G-E_{\overline{I_Y}\cap
\overline{\{j\}}}-B_{\overline{\{j\}}-\overline{I_Y}}-
A_{\overline{I_Y}-\overline{\{j\}}}).\eqno{(11)}$$

It is not difficult to prove the following two claims:\\

{\bf Claim 1}
$$\{j\}\triangle I_Y=\overline{\{j\}}\triangle \overline{I_Y}.$$

{\bf Claim 2}\ \ The mapping $\phi: \{b_i'|\ i\in
I_{Y'}\}\longmapsto \{a_i|\ i\in I_{Y'}\}$ is a bijection between
the set of the odd subsets of $B'$ and the set of the odd subsets
of $A$.\\

By Claims $1-2$ and $(10)-(11)$, the following claim is
obvious:\\

{\bf Claim 3} \[\sum_{{Y'\subseteq B'}\atop{|Y'|\ \mbox{is odd}}}
M(G'-a_j'-Y') M(G'-A'\backslash \{a_j'\}-\overline {Y'})=
\sum_{{Y\subseteq B}\atop {|Y|\ \mbox{\small is
odd}}}\left(\prod_{e\in E_{\{j\}\triangle
I_Y}}\omega_e\right)\times \]\[\left\{
M(G-E_{I_Y\cap\{j\}}-B_{\{j\}-I_Y}-A_{I_Y-\{j\}}) M(G-E_{\overline
{I_Y}\cap \overline {\{j\}}})-B_{\overline {\{j\}}-\overline
{I_Y}}-A_{\overline {I_Y}-\overline {\{j\}}})\right\}.\]\\

\par Let $W'$ be an even subset of $B'$. By our notation defined
above, $I_{W'}=\{i|\ b_i'\in W'\}$. Let $W=\{b_i|\ b_i'\in W'\}$,
$I_W=I_{W'}$. As in the proof of Claim 3 we can prove the
following claim:\\

{\bf Claim 4} \[\sum_{{W'\subseteq B'}\atop {|W|\ \mbox{is
even}}}M(G'-W')M(G'-A'-\overline {W'})=\sum_{{W\subseteq B}\atop
{|W|\ \mbox{\small is even}}}\left(\prod_{e\in
E_{I_W}}\omega_e\right)M(G-A_{I_W})M(G-E_{\overline{I_W}}-B_{I_W})\]

The theorem is immediate from Claims $3-4$ and $(9)$. \ \ \ \
$\blacksquare$\\

\par If we set $k=2$, it is not difficult to see that the following corollary holds.\\

{\bf Corollary 3.3}\ \ Let $G$ be a plane weighted graph with even
number of vertices. Let $e_1=a_1b_1$ and $e_2=a_2b_2$ be two
independent edges on the boundary of a face $f$ of $G$ and $a_1,
b_1, a_2, b_2$ appear in a cyclic order on a face of $G$. Then
$$M(G)M(G-e_1-e_2)+\omega_{e_1}\omega_{e_2}M(G-a_1-a_2)M(G-b_1-b_2)$$
$$=M(G-e_1)M(G-e_2)+\omega_{e_1}\omega_{e_2}M(G-a_1-b_2)M(G-a_2-b_1)
,$$
where $\omega_{e}$ denotes the weight of edge $e$.\\

\par {\bf Corollary 3.4}\ \ Let $G=(U,V)$ be a plane weighted bipartite
graph, in which $|U|=|V|=n$ and the weight of every edge $e$ in
$G$ is denoted by $\omega_e$. Let $e_1=a_1b_1, e_2=a_2b_2, \ldots,
e_k=a_kb_k$ ($2\leq k\leq n$) be k independent edges in the
boundary of a face $f$ of $G$ and let vertices $a_1, b_1, a_2,
b_2, \ldots, a_k, b_k$ appear in a cyclic order on $f$. If
$A=\{a_i|1\leq i\leq k\}\subseteq U$ and $B=\{b_i|1\leq i\leq
k\}\subseteq V$, then for any $j=1,2,\ldots,k$\\

$$M(G)M(G-e_1-e_2-\cdots-e_k)=$$
\[M(G-e_j)M(G-\overline{\{e_j\}})+\sum_{{1\le i\le k}\atop {i\neq j}}
\omega_{e_i}\omega_{e_j}M(G-a_i-b_j)M(G-a_j-b_i-E_{\overline{\{i,j\}}}),\
\ (12)\] where
$\overline{\{e_j\}}=\{e_1,e_2,\ldots,e_k\}\backslash \{e_j\}$ and
$E_{\overline{\{i,j\}}}=\{e_t|t\in {\bf[k]}\backslash
\{i,j\}\}$.\\

{\bf Proof}\ \ Note that if $W$ is a nonempty even subset of $B$
or $Y$ is an odd subset of $A$ such that $|Y|\geq 3$ then
$M(G-A_{I_W})=0$ and
$M(G-E_{I_Y\cap\{j\}}-B_{\{j\}-I_Y}-A_{I_Y-\{j\}})=0$ in $(8)$ in
Theorem 3.2 (since $G$ is a bipartite graph, and $A\subseteq U,
B\subseteq V$). Hence the corollary is immediate from Theorem
3.2.\ \ \ \ $\blacksquare$\\

\par One direct corollary of Corollaries 3.4 is the following
result:\\

{\bf Corollary 3.5}\ \ Let $G=(U,V)$ be a plane weighted bipartite
graph in which $|U|=|V|$. Let $e_1=a_1b_1$ and $e_2=a_2b_2$ be two
independent edges on the boundary of a face $f$ of $G$ and
$a_1, b_1, a_2, b_2$ appear in a cyclic order on a face of $G$. \\

(1)\ \ If $\{a_1,a_2\}\subseteq U$ and $\{b_1,b_2\}\subseteq V$,
then
\begin{center}
$M(G)M(G-e_1-e_2)=M(G-e_1)M(G-e_2)+\omega_{e_1}\omega_{e_2}M(G-a_1-b_2)M(G-a_2-b_1).$
\end{center}

(2)\ \ If $a_1\in U$ and $a_2\in V$ or $a_1\in V$ and $a_2\in U$,
then
$$M(G)M(G-e_1-e_2)=M(G-e_1)M(G-e_2)-\omega_{e_1}\omega_{e_2}
M(G-a_1-a_2)M(G-b_1-b_2),$$
where $\omega_{e}$ denotes the weight of edge $e$.\\

\par By the method similar to that in the proof of Theorem 3.2, we can prove the
following result:\\

{\bf Theorem 3.6}\ \ Let $G$ be a plane weighted graph with even
number of vertices. Let $a_1$ and $b_1$ be two vertices of $G$ and
$e=a_2b_2$ an edge of $G$. If
the four vertices $a_1,b_1,a_2,b_2$ appear in a cyclic order on a face of $G$, then \\
$$M(G)M(G-a_1-b_1-e)=$$
$$
M(G-a_1-b_1)M(G-e)+\omega_eM(G-a_1-a_2)M(G-b_1-b_2)-\omega_eM(G-a_1-b_2)
M(G-a_2-b_1).
$$
\par A direct corollary of Theorem 3.6 is the following result:\\

{\bf Corollary 3.7}\ \ Let $G=(U,V)$ be a plane weighted bipartite
graph in which $|U|=|V|$. Let $a_1$ and $b_1$ be two vertices of
$G$ with different colors and $e=a_2b_2$ an edge of $G$. If
$a_1,b_1,a_2,b_2$ appear in a cyclic order of
a face of $G$, then \\
$(i)$\ \ if $\{a_1,b_2\}\subseteq U$ and $\{a_2,b_1\}\subseteq V$
(or $\{a_1,b_2\}\subseteq V$ and $\{a_2,b_1\}\subseteq U$) then
$$M(G)M(G-a_1-b_1-e)=M(G-a_1-b_1)M(G-e)+\omega_eM(G-a_1-a_2)M(G-b_1-b_2);$$
$(ii)$\ \ if $\{a_1,a_2\}\subseteq U$ and $\{b_1,b_2\}\subseteq V$
or $\{a_1,a_2\}\subseteq V$ and $\{b_1,b_2\}\subseteq U$ then
$$M(G)M(G-a_1-b_1-e)=M(G-a_1-b_1)M(G-e)-\omega_eM(G-a_2-b_1)M(G-a_1-b_2);$$
where $\omega_e$ is the weight of edge $e=a_2b_2$.\\

{\bf Remark 4}\ \ Let $G=(U,V)$ be a plane weighted graph with
even number of vertices. Let $a_i$ and $b_i$ for $i=1,2,\ldots,s$
be $2s$ vertices of $G$, and let $e_i=a_{s+i}b_{s+i}$ for
$i=1,2,\ldots,t$ be $t$ edges of $G$ ($6\leq s+t\leq n$). If
vertices $a_1,b_1,a_2,b_2,\ldots,a_{s+t},b_{s+t}$ appear in the
boundary of a face $f$ of $G$(which may appear in different order
of $f$), we can consider the problems similar to Theorem 3.6.

\section* { \bf 4. Weighted Aztec diamonds}
\ \ \ \ In this section, we use Corollary 3.5 to give a new proof
of one identity concerning perfect matchings of the weighted Aztec
diamond in Yan and Zhang \cite{YZ05}, which implies a formula on
the sum of weights of perfect matchings of the weighted Aztec
diamond in \cite{Ciucu98,Stanley}.\\

The Aztec diamond of order $n$, denoted $AD_n$, is defined to be
the graph whose vertices are the white squares of a $(2n+1)\times
(2n+1)$ chessboard with black corners, and whose edges connect
precisely those pairs of white squares that are diagonally
adjacent (Figure 5(a) illustrates $AD_4$). In \cite{EKLP92}, four
proofs are presented that $M(AD_n)=2^{n(n+1)/2}$. Ciucu
\cite{Ciucu98} showed that $M(AD_n) = 2^n M(AD_{n-1})$, which
clearly implies the previous formula (since $M(AD_1)=2$). By two
different methods, Kuo \cite{Kuo04} and Yan and Zhang \cite{YZ05}
proved that
$$M(AD_n)=\frac{2 M(AD_{n-1})^2}{M(AD_{n-2})}\eqno{(13)}$$
which, in turn, implies that $M(AD_n)=2^{n(n+1)/2}$. Recently, Eu,
Fu \cite{EF05} and Brualdi and Kirkland \cite{BK05} gave
independently a new method to prove this formula.\\

\par Stanley weighted the Aztec diamond of order $n$ as follows. Weight every
$4-$cycle in the $i$th column by assigning the variables
$x_i,y_i,w_i$ and $z_i$ to its four edges, starting with the
northwestern edge and going clockwise. We denote this weighted
Aztec diamond of order $n$ by $(AD_n;1\leq i\leq n)$. The case
$n=4$, i.e. $(AD_4,1\leq i\leq 4)$, is illustrated in Figure 5(a),
and the array on the right indicates the weight pattern on edges.
We can also weight every $4-$cycle of $AD_n$ in the $i$th column
by assigning the variables $x_{i+1},y_{i+1},w_{i+1}$ and $z_{i+1}$
to its four edges, starting with northwestern edge and going
clockwise. Denote this weight Aztec diamond of order $n$ by
$(AD_n;2\leq i\leq n+1)$. The case n=3, i.e. $(AD_3,2\leq i\leq
4)$, is illustrated in Figure 5(b), and the array on the right
indicates the weight pattern on the edges).\\

Based on the method on the graphical vertex-condensation Yan and
Zhang \cite{YZ05} proved that
$$M(AD_n;1\leq i\leq n)M(AD_{n-2};2\leq i\leq n-1)=$$
$$(x_1w_{n}+y_nz_1)
M(AD_{n-1};1\leq i\leq n-1)M(AD_{n-1};2\leq i\leq n),\eqno{(14)}$$
which implies the following theorem by induction on $n$, which was
previously proved by Stanley \cite{Stanley} and Ciucu
\cite{Ciucu98}.
%%%%%%%%%%%%%%%%%%%%%%%%%%%%%%%%%%%%%%%%%%
%%%%%%%%%%%%% Figure 5
%%%%%%%%%%%%%%%%%%%%%%%%%%%%%%%%%%%%%%%%%%
\begin{figure}[htbp]
  \centering
  \scalebox{0.8}{\includegraphics{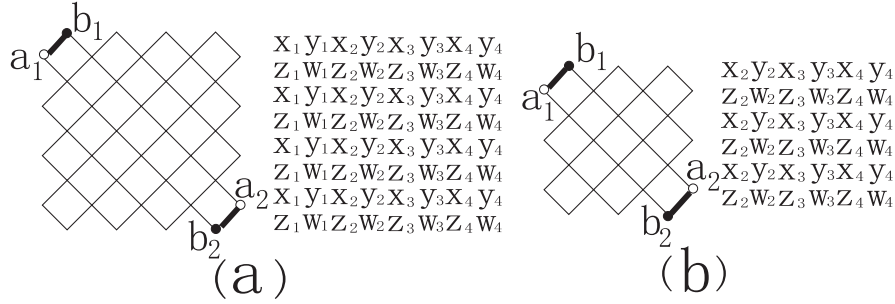}}
  \caption{\ (a)\ The weighted Aztec diamond
$(AD_4;1\leq i\leq 4)$. \ (b)\ The weighted Aztec diamond
$(AD_3;2\leq i\leq 4)$.}
\end{figure}

{\bf Theorem 4.1} (Stanley \cite{Stanley} and Ciucu
\cite{Ciucu98})\ \ The sum of weights of perfect matchings of the
weighted Aztec diamond $(AD_n;1\leq i\leq n)$ of order $n$
$$M(AD_n;1\leq i\leq n)=\prod_{1\leq i\leq j\leq n}(x_iw_j+z_iy_j).$$

Now we use Corollary 3.5 to give a new proof of $(14)$ as
follows.\\

Let $G=(AD_n;1\leq i\leq n)$. For the sake of convenience, we
rotate clockwise $AD_n$ by $45^\circ$ so that their edges are
horizontal and vertical. Let $a_1$ and $b_1$ be the two vertices
which are the left and right vertices of the horizontal edge in
the northern corner, and let $a_2$ and $b_2$ be the two vertices
which are the right and left vertices of the horizontal edge in
the southern corner respectively. The cases $n=3$ and $4$ rotated
by $45^\circ$ are illustrated in Figure 5(b) and (a),
respectively. Obviously, two edges $e_1=a_1b_1$ and $e_2=a_2b_2$
appear the boundary of the unbounded face of $G$. Particularly,
$a_1$ and $a_2$ share one color, and $b_1$ and $b_2$ have another
color. Then, by Corollary 3.5, we have
$$M(G)M(G-e_1-e_2)=M(G-e_1)M(G-e_2)+\omega_{e_1}\omega_{e_2}M(G-a_1-b_2)M(G-a_2-b_1).
\eqno {(15)}$$
Note that, after the removing the forced edges, we
have
$$
M(G-e_1-e_2)=(y_nz_1)^{n-1}(y_1y_2\ldots y_n)(z_1z_2\ldots
z_n)M(AD_{n-2}; 2\leq i\leq n-1);\eqno{(16)}
$$
$$
M(G-e_1)=z_1^{n}(y_1y_2\ldots y_n)M(AD_{n-1};2\leq i\leq n);
\eqno{(17)}
$$
$$
M(G-e_2)=y_n^n(z_1z_2\ldots z_n)M(AD_{n-1}; 1\leq i\leq n-1);
\eqno{(18)}
$$
$$
M(G-a_1-b_2)=y_n^{n-1}(y_1y_2\ldots y_n)M(AD_{n-1};2\leq i\leq n);
\eqno{(19)}
$$
$$
M(G-a_2-b_1)=z_1^{n-1}(z_1z_2\ldots z_n)M(AD_{n-1};1\leq i\leq
n-1). \eqno{(20)}
$$
Note that $\omega_{e_1}=x_1$ and $\omega_{e_1}=w_n$. Hence $(14)$
is immediate from $(15)-(20)$.

\vskip 0.3cm \noindent {\bf Acknowledgements} \vskip 0.5cm
\noindent

%\par This work is completed during the first author's stay
%in Institute of Mathematics,Academia Sinica, Taipei. The first author
%thank the Institute for its support.
Thanks to all people (such as James Propp, Rick Kenyon, Eric
Heng-Shiang Kuo, Christian Krattenthaler, etc.) for the full
discussion from the ``domino archives". Particularly, thanks to
James Propp for some helpful suggestions for this paper. Thanks
also to the referees for providing some very helpful suggestions
for this paper. Professor Krattenthaler have told us by an E-mail
that he could use the Pfaffian method to prove some identities as
in Theorem 2.2. \vskip1cm \noindent
\newcounter{cankao}
\begin{list}
{[\arabic{cankao}]}{\usecounter{cankao}\itemsep=0cm}
\centerline{\bf References} \vspace*{0.5cm} \small
%\bibitem{BE04}
%G. Benkart and O. Eng, Weighted Aztec Diamond Graphs and the Weyl
%Character Formula, Electronic Journal of Combinatorics, 11(2004),
%\#R28.
\bibitem{BK05}
R. Brualdi and S. Kirkland, Aztec diamonds and digraphs, and
Hankel determinants of Schr\"oder numbers, J. Combin. Theory Ser.
B 94(2005), 334$-$351.
\bibitem{Ciucu97}
M. Ciucu, Enumeration of Perfect Matchings in Graphs with
Reflective Symmetry, J. Combin. Theory Ser. A, 77(1997), 67$-$97.
\bibitem{Ciucu96}
M. Ciucu, Enumeration of perfect matchings of cellular graphs, J.
Algebraic Combin. 5(1996), 87$-$103.
\bibitem{Ciucu98}
M. Ciucu, A complementation theorem for perfect matchings of
graphs having a cellular completion, J. Combin. Theory Ser. A,
81(1998), 34$-$68.
\bibitem{DW95}
A. W. M. Dress and W. Wenzel, A simple proof of an identity
concerning Pfaffians of skew symmetric matrices, Adv. Math.,
112(1995), 120$-$134.
\bibitem{EKLP92}
N. Elkies, G. Kuperberg, M. Larsen, and J. Propp,
Alternating$-$sign matrices and domino tilings (Parts I and II),
J. Algebraic Combin. 1(1992), 111$-$132 and 219$-$234.
\bibitem{EF05}
S. P. Eu and T. S. Fu, A simple proof of the Aztec diamond
theorem, Electron. J. Combin. 12(2005), \#R18.
\bibitem{Ham01}
A. M. Hamel, Pfaffian identities: a combinatorial approach, J.
Combin. Theory, Ser. A, 94(2001), 205$-$217.
\bibitem{Kas63}
P.W.Kasteleyn, Dimer statistics and phase transition, J. Math.
Phys. 4 (1963), 287-293.
\bibitem{Kas67}
P.W.Kasteleyn, Graph Theory and Crystal Physics, Graph Theory and
Theoretical Physics (F.Harary, ed.), Academic Press, 1967, 43-110.
\bibitem{Knuth96}
D. E. Knuth, Overlapping pfaffians, Electron. J. Combin., 3(1996),
{\bf R}5.
\bibitem{Kuo04}
E. H. Kuo, Applications of Graphical Condensation for Enumerating
Matchings and Tilings, Theoret. Comput. Sci., 319 (2004), 29-57.
%\bibitem{Mac78}
%P. MacMahon, Menoir on the Theory of Partitions of Numbers$-$Part
%V. Partitions in Two$-$dimension Space, Philosophical Transactions
%of the Royal Society of London, 211(1912), 75$-$110. Reprinted in
%Percy Alexander MacMahon: Collected Papers, ed. George E. Andrews,
%Vol. 1, pp. 1328$-$1363. MIT Press, Cambridge, Mass., 1978.
\bibitem{Propp03}
J. Propp, Generalized Domino$-$Shuffling, Theoret. Comput. Sci.,
303(2003), 267$-$301.
%\bibitem{Propp97}
%J. Propp, Talk, American Mathematical Society Meeting, San Diego,
%CA, Jan. 1997.
\bibitem{Stanley}
R. P. Stanley, Private communication.
\bibitem{Wenzel93}
W. Wenzel, Pfaffian forms and $\Delta-$matroids, Discrete Math.,
115(1993), 253$-$266.
\bibitem{YZ05}
W. G. Yan and F. J. Zhang, Graphical Condensation for Enumerating
Perfect Matchings, J. Combin. Theory Ser. A, 110(2005), 113-125.
\end{list}
\end{document}